\documentclass[11pt, twoside]{article}

\usepackage{amsmath,amssymb,url}
\usepackage{enumerate, graphicx, array}

\usepackage{a4wide,fancyhdr}
\newtheorem{theorem}{Theorem}[section]
\newtheorem{lemma}[theorem]{Lemma}
\newtheorem{corollary}[theorem]{Corollary}
\newtheorem{proposition}[theorem]{Proposition}

\newtheorem{definition}[theorem]{Definition}
\newtheorem{remarks}[theorem]{Remarks}

\def\bproof{{\noindent\bf Proof: }}
\def\eproof{\hfill$\Box$}
\newcommand\into{\hookrightarrow}
\newcommand\onto{\twoheadrightarrow}
\newcommand\mapsfrom{\mathrel{\reflectbox{\ensuremath{\mapsto}}}} %%needs graphicx

\DeclareMathOperator{\Con}{Con}

\begin{document}

\title{Absolute Retracts and Essential Extensions in Congruence Modular Varieties}

\author{Peter Ouwehand}
\date{15 December 2011}
 \maketitle

\begin{abstract} This paper studies absolute retracts in congruence modular varieties of universal algebras. It is shown that every absolute retract with finite dimensional congruence lattice is a product of subdirectly irreducible algebras. Further, every absolute retract in a residually small variety is the product of an abelian algebra and a centerless algebra.
\end{abstract}

\fancyhead[LE,RO]{\thepage} \fancyhead[LO]{P. Ouwehand}\fancyhead[RE]{\em Absolute Retracts in Modular Varieties}

\section{Introduction}
Recall that an algebra $R$ belonging to a variety $\mathcal V$ of algebras is said to be an {\em absolute retract} in $\mathcal V$ if and only if it is a retract of each of its extensions in $\mathcal V$, that is, if for any embedding $e\colon R\into A$ in $\mathcal V$ there is a surjective homomorphism $p\colon A\onto R$ such that $p\circ e = \text{id}_R$. 
Absolute retracts have been studied in a universal-algebraic context in a number of papers, particularly in association with related notions of injectivity, congruence extension, amalgamation \cite{GL:1972}, \cite{Bergman:1985}, \cite{DK:1986}, \cite{JR:1989}, \cite{OR:1996}, \cite{JJOR:2001}, \cite{Ouwehand:2009}. 

We briefly describe why the investigation of absolute retracts is worth pursuing: In his survey on equational logic, Walter Taylor\cite{Taylor:1979} states the following problem: {\em Does every residually small variety of algebras with the amalgamation property also have the congruence extension property?} This problem is still open, though it has been settled in the affirmative for congruence modular varieties by the efforts of  Bergman, Kearnes and McKenzie \cite{Bergman:1986}, \cite{BM:1988}, \cite{Kearnes:1989}. Absolute retracts seem central to this problem. Firstly, the notion of absolute retract joins two model-theoretic properties: An algebra is an absolute retract in a variety if and only if it is both equationally compact and algebraically closed in that variety (cf. \cite{Ouwehand:2009}). Secondly, every absolute retract in a variety is a member of the amalgamation class of that variety (cf. \cite{Bergman:1985}). Thirdly, a variety $\mathcal V$ is residually small if  and only if every algebra in $\mathcal V$ is embeddable into an absolute retract of $\mathcal V$ (cf. \cite{Taylor:1972}). Finally, it is easy to show that if a residually small variety $\mathcal V$ satisfies the amalgamation property, then $\mathcal V$ satisfies the congruence extension property if and only if the class of absolute retracts of $\mathcal V$ is closed under products.

The investigation of absolute retracts has been most fruitful when restricted to congruence distributive varieties, where the Fraser-Horn property and J\'onsson's Lemma provide powerful tools for managing congruences on products.  Davey and Werner\cite{DW:1979} observe that ``there are a number of papers in which it is shown that the injectives, and more generally the weak injectives, of a particular variety are precisely finite complete Boolean powers of appropriate subdirectly irreducible algebras", and proceed to prove a statement which encompasses many such results at a single stroke.
Products of Boolean powers also play a role in \cite{OR:1996} and \cite{JJOR:2001}, where it is proved that, in a finitely generated congruence distributive variety, every absolute retract is a product of Boolean powers (indeed, reduced powers) of maximal subdirectly irreducible algebras.

In the congruence modular case, congruences on products are less manageable. Davey and Kov\'acs \cite{DK:1986} {\em do} study absolute retracts in modular varieties, but restrict their attention to directly indecomposable
absolute retracts, where they show that each such is finitely subdirectly irreducible. The main aim of the current paper is to partly extend this result to cover products: We prove that if $R$ is an absolute retract in a congruence modular variety, and if the congruence lattice of $R$ is of finite dimensional, then $R$ is a product of subdirectly irreducible algebras.

\section{Preliminaries}
In this section we fix notation and state some basic results that we require in the sequel.

We denote the congruence lattice of an algebra $A$ by $\Con(A)$. $0_A$ and $1_A$ denote the smallest and largest congruences on $A$, but we may omit the subscripts if these are clear from context.  If $\varphi\leq\psi\in\Con(A)$, then $I[\varphi,\psi]$ denotes the interval $\{\theta\in\Con(A):\varphi\leq \theta\leq \psi\}$ in $\Con(A)$. If $a\in A$ and $\theta\in\Con(A)$ we shall use two notations for the congruence class of $a$ modulo $\theta$: We shall denote by $[a]\theta$ the congruence class of $a$ as a {\em set}: $[a]\theta:=\{b\in A: a\,\theta\,b\}$. On the other hand, the {\em element} of the quotient algebra $A/\theta$ corresponding to $a$ will be denoted by $a/\theta$. If $B\subseteq A$, then we define $[B]\theta:=\bigcup_{b\in B}[b]\theta$. 
If $A$ is a subalgebra of $B$ and $\theta\in\Con(B)$, then we shall denote the restriction of $\theta$ to $A$ by $\theta\restriction_A$.

If $\theta_i\in\Con(A_i)$ for $i\in I$, then  $\prod_{i\in I}\theta_i$ is the congruence on $\prod_IA_i$ given by: $(a_i)_{i\in I}\,(\prod_{i\in I}\theta_i)\, (b_i)_{i\in I}$ if and only if $a_i\,\theta_i\,b_i$ for all $i\in I$. For finite products, we may write $\theta_1\times\theta_2\times\dots\times\theta_n$ instead of $\prod_{i=1}^n\theta_i$. A congruence of the form $\prod_{i\in I}\theta_i$ is called a {\em product congruence}.

The following result is implicit in \cite{DK:1986}:
\begin{lemma}\label{lemma_product_congruence} Suppose that $\Con(A\times B)$ is a modular lattice, and that $\varphi_1\leq \varphi_2\in\Con(A)$ and $\psi\in\Con(B)$. Then each $\theta\in\Con(A\times B)$ satisfying $\varphi_1\times\psi\leq\theta\leq \varphi_2\times\psi$ is a product congruence of the form $\varphi\times\psi$ for some $\varphi_1\leq\varphi\leq\varphi_2$.
\end{lemma}

\bproof By modularity, since $(\varphi_2\times\psi)\land (\varphi_1\times1_B) =\varphi_1\times\psi$ and $(\varphi_2\times\psi)\lor (\varphi_1\times1_B) =\varphi_2\times1_B$, we have $I[\varphi_1\times\psi,\varphi_2\times\psi]\cong I[\varphi_1\times 1_B, \varphi_2\times 1_B]$, via inverse isomorphisms
\[\theta\mapsto \theta\lor(\varphi_1\times 1_B)\qquad
\gamma\land(\varphi_2\times\psi)\mapsfrom \gamma\] Obviously, $I[\varphi_1\times 1_B, \varphi_2\times 1_B]\cong I[\varphi_1,\varphi_2]$. Thus whenever $\theta\in I[\varphi_1\times \psi,\varphi_2,\times\psi]$, then $\theta\lor (\varphi_1\times 1_B) =\varphi\times 1_B$ for some $\varphi\in I[\varphi_1,\varphi_2]$, and hence $\theta =(\varphi\times 1_B)\land (\varphi_2\times\psi) = \varphi\times\psi$.
\eproof

An algebra $A$ is said to be {\em finitely subdirectly irreducible} if and only if $0_A$ is meet irreducible. 
We say that a congruence $\alpha\in\Con(A)$ is {\em dense} if and only if whenever $\theta\not=0_A$ in $\Con(A)$, then  $\theta\land\alpha\not=0_A$ also. Clearly, if $A$ is finitely subdirectly irreducible, then every non-zero congruence on $A$‎ is dense. 

\begin{lemma}\label{lemma_prod_dense} Let $A, B$ be algebras such that $\Con(A\times B)$ is a modular lattice.  If $\alpha$ is dense in $\Con(A)$, and $\beta$ is dense in $\Con(B)$, then $\alpha\times\beta$ is dense in $\Con(A\times B)$.\end{lemma}
\bproof We first show that $\alpha\times 1$ is dense in $\Con(A\times B)$: For if $\theta\land(\alpha\times 1) = 0$, then $[(\alpha\times 1))\land\theta]\lor(0\times 1) =(0\times 1)$, and hence $(\alpha\times1)\land[\theta\lor(0\times 1)]=0\times 1$, by modularity. But as $\theta\lor(0\times 1)\in I[0\times 1, 1\times 1]$, there is by Lemma \ref{lemma_product_congruence}, an $\bar{\alpha}\in\Con(A)$ such that $\theta\lor(0\times 1) =\bar\alpha\times 1$.
Thus $(\alpha\land\bar\alpha)\times 1=0\times 1$, 
so that we conclude first that $\bar\alpha\land\alpha = 0$, and then  that $\bar\alpha = 0$, because $\alpha$ is dense. Then $\theta\lor(0\times 1)= 0\times 1$, and hence $\theta\leq0\times 1\leq\alpha\times 1$. Thus $0=\theta\land(\alpha\times 1) =\theta$, as required.

In the same way it follows that $1\times\beta$ is dense in $\Con(A\times B)$. Now observe that $(\alpha\times \beta) = (\alpha\times 1)\land(1\times\beta)$, and conclude that $\alpha\times \beta$ is dense in $\Con(A\times B)$.
\eproof

An extension $A\into B$ is said to be an {\em essential extension} if and only if whenever $\theta\in\Con(B)$ has $\theta\restriction_A = 0_A$, then $\theta=0_B$. The following facts are well--known (and easy to show)(cf. \cite{Bergman:1985}): \begin{enumerate}[(i)]\item An algebra $R\in\mathcal V$ is an absolute retract in $\mathcal V$  if and only if it has no proper essential extensions in $\mathcal V$. \item An essential extension of a (finitely) subdirectly irreducible algebra is also (finitely) subdirectly irreducible.\end{enumerate} A subdirectly irreducible algebra is said to be a {\em maximal subdirectly irreducible} in $\mathcal V$ if and only if it has no proper essential extensions in $\mathcal V$. Thus each maximal subdirectly irreducible is an absolute retract.

Henceforth, we work in a congruence modular variety $\mathcal V$, where we shall make use of the commutator theory. We enumerate the following basic facts, which are proved in \cite{FM:1987}, \cite{Gumm:1983}:
\begin{enumerate}[1.]\item\label{fact_basic_comm} If $A\in\mathcal V$ and $\varphi,\psi,\psi_i\;(i\in I)\in\Con(A)$, then
\begin{enumerate}[(i)]\item $[\varphi,\psi]\leq\varphi\land\psi$;
\item $[\varphi,\psi] = [\psi,\varphi]$;
\item $[\varphi,\bigvee_{i\in I}\psi_i] = \bigvee_{i\in I}[\varphi,\psi_i]$;
\end{enumerate}\item \label{fact_comm_restr} If $B$ is a subalgebra of $A$ and $\varphi,\psi\in\Con(A)$, then $[\varphi\restriction_B,\psi\restriction_B]\leq [\varphi,\psi]\restriction_B$.
\item \label{fact_comm_image} If $\pi\in\Con(A)$, then, with $\varphi,\psi\geq\pi$, we have that $\left[\varphi/\pi,\psi/\pi\right] =\dfrac{[\varphi,\psi]\lor\pi}{\pi}$ in $\Con(A/\pi)$.
\item \label{facty_comm_center}
A congruence $\theta\in\Con(A)$ is said to be {\em abelian} if $[\theta,\theta]=0_A$, and {\em central} if $[\theta,1_A]=0_A$. Clearly every central congruence is abelian. The {\em center} of $A$ is the largest central congruence of $A$, and denoted by $\zeta_A$. The algebra $A$ is said to be {\em centerless} if $\zeta_A = 0_A$.
Furthermore:
\begin{enumerate}[(i)]\item $\zeta_{A\times B} =\zeta_A\times\zeta_B$;
\item If $B$ is a subalgebra of $A$, then $\zeta_A\restriction_B\leq \zeta_B$. In particular, the restriction of a central congruence is itself central.
\item Abelian congruences permute with every congruence: If $\theta,\varphi\in\Con(A)$ and $\theta$ is abelian, then $\theta\circ\varphi=\varphi\circ\theta=\theta\lor\varphi$.
\end{enumerate}
\item\label{fact_comm_diff}
Each congruence modular variety $\mathcal V$ has a ternary term $d(x,y,z)$, called the {\em  Gumm difference term}, with the following properties:
\begin{enumerate}[(i)]\item $\mathcal V\vDash d(x,y,y) = x$;
\item  If $\theta$ is abelian, and $x\,\theta\,y$, then  $d(x,x,y) = y$.
\end{enumerate}
\item \label{fact_comm_C1} A residually small congruence modular variety satisfies the commutator equation (C1) (cf. \cite{FM:1987}, Theorems 8.1 and 10.4):
\[\alpha\land[\beta,\beta] = [\alpha\land\beta, \beta]\tag{C1}\]
\end{enumerate}

We will also need the following result:

\begin{theorem} \label{thm_Gumm_9.1}{\rm (\cite{Gumm:1983}, Thm. 9.1)} Suppose that $\mathcal V$ is a congruence modular variety with Gumm difference term $d(x,y,z)$. Let $A\in\mathcal V$, and $\varphi\geq\psi\in\Con(A)$. Then $[\varphi,\psi]=0_A$ if and only if for every term operation $f$ (of arity $n$) and every $x_i\,\varphi\, y_i\,\psi \,z_i$ (for $1\leq i\leq n$) we have
\begin{enumerate}[(i)]\item $d(y_i,y_i,z_i) = z_i$, and
\item \label{d_commutes}$f\left(d\left(\begin{matrix}x_1\\y_1\\z_1\end{matrix}\right), d\left(\begin{matrix}x_2\\y_2\\z_2\end{matrix}\right),\dots, d\left(\begin{matrix}x_n\\y_n\\z_n\end{matrix}\right)\right) = d\left(\begin{matrix}f(x_1,x_2,\dots, x_n)\\f(y_1,y_2,\dots, y_n)\\f(z_1,z_2,\dots,z_n)\end{matrix}\right)$
\end{enumerate}
\end{theorem}

\begin{corollary}\label{corr_p_id} If $x\,\zeta_A\,y$, and $a\in A$, then 
$d(x,a,d(a,x,y)) = y$.
\end{corollary}
\bproof
Note that $1_A\geq\zeta_A$, and that $[1,\zeta_A]=0$. As $x\,\zeta_A\,y$, we may apply Theorem \ref{thm_Gumm_9.1} with $f=d$ to conclude that $d$ commutes with itself on $\left(\begin{matrix}x&x&x\\a&x&x\\a&x&y\end{matrix}\right)$, i.e.  that $d\left(\begin{matrix} d(x,x,x)\\d(a,x,x)\\d(a,x,y)\end{matrix}\right)
= d\left(\begin{matrix} d\left(\begin{matrix}x\\a\\a\end{matrix}\right), d\left(\begin{matrix}x\\x\\x\end{matrix}\right),
d\left(\begin{matrix}x\\x\\y\end{matrix}\right)\end{matrix}\right)$, so that $d(x,a,d(x,a,y))=d(x,x,d(x,x,y)$.
But as $\zeta_A$ is abelian, and $x\zeta_Ay$, we have $d(x,x,y)= y$, and hence $d(x,a,d(a,x,y)) = y$.
\eproof

\section{Subdirect Product--Essential Extensions are Essential}
We work throughout in a congruence modular variety $\mathcal V$ with Gumm difference term $d$.

\begin{definition} \rm Let $n\in\mathbb N$.  An embedding $e\colon A\into\prod_{i=1}^nA_i$ is said to be a {\em product--essential extension} if and only if whenever $\varphi_i\in\Con(A_i)$ are such that $(\prod_{i=1}^n\varphi_i)\restriction_A=0_A$, then $\varphi_i=0_{A_i}$ (for $i=1,2,\dots,n$).
%\endbox
\end{definition}

\begin{remarks}\label{rem_product_essential}\rm
\begin{enumerate}[(a)]\item\label{rem_eta_max}
Note that if $e\colon A\into \prod_{i=1}^nA_i$ is a subdirect product\--essential extension, and if $\eta_i:=\ker \pi_i\restriction_A$ is the kernel of the natural projection $\pi_i\circ e\colon A\onto A_i$, then 
\begin{enumerate}[(i)]\item $\bigwedge_{i=1}^n\eta_i=0_A$;
\item If $\eta_i\leq\phi_i$ for $i=1,\dots, n$, and $\bigwedge_{i=1}^n\phi_i=0_A$, then $\phi_i=\eta_i$ for $i=1,\dots, n$.
\end{enumerate}
\item \label{rem_sd_pr_ess} If $A\into\prod_{i=1}^nA_i$ is a subdirect embedding, then there are $\bar\varphi_i\in\Con(A_i)$ such that $A\into \prod_{i=1}^nA_i/\bar\varphi_i$ is a subdirect product--essential embedding. Indeed, if $A_i\cong A/\eta_i$, then by Zorn's lemma there is  $\varphi_1\in\Con(A)$ maximal with respect to the properties that $\varphi_1\geq\eta_1$ {\em and} $\varphi_1\land\eta_2\dots\land \eta_n = 0$.  Then choose $\varphi_2\in\Con(A)$ maximal such that $\varphi_2\geq \eta_2$ and $\varphi_1\land\varphi_2\land\eta_3\dots
\land \eta_n = 0$, etc.  Then $\bar\varphi_i:=\varphi_i/\eta_i$ will do.
\end{enumerate}%\endbox
\end{remarks}

The main technical result of this section is the following:
\begin{theorem} \label{thm_sd_pr_ess} If $e\colon A\into\prod_{i=1}^n A_i$ is a subdirect product--essential embedding in a congruence modular variety, then it is an essential embedding.
\end{theorem}

Before we can tackle the proof, we  will require some intermediary results, many of which  owe a large debt to Davey and Kov\'acs\cite{DK:1986}.

\begin{proposition}\label{propn_central_essential}  Let $A$ be a subalgebra of an algebra $B$, and suppose that $\theta\in\Con(B)$ is a central congruence in $B$ with $\theta\restriction_A=0_A$ and $[A]\theta=A$. Then $\theta=0_B$.
\end{proposition}

\bproof Let $a\in A$. If $x\,\theta\, y$, then $d(a,x,y)\,\theta\, d(a,x,x) = a$, and hence $d(a,x,y)\in A$. But then $(d(a,x,y),a)\in\theta\restriction_A=0_A$, so $d(a,x,y) = a$. Thus $x=d(x,a,a) = d(x,a,d(a,x,y)) = y$, by Corollary \ref{corr_p_id}.
\eproof

\begin{proposition}\label{propn_theta_central} Suppose that $f\colon B\into\prod_{i=1}^nB_i$ is a product--essential embedding. Further suppose that $\bar\theta\in\Con(\prod_{i=1}^n B_i)$ is such that $\bar\theta\restriction_B=0_B$.  Then $\bar\theta$ is a central congruence.
\end{proposition} 

\bproof Define $\bar\nu_i=0\times\dots\times0\times1\times0\times\dots0$, where $1$ occurs in the $i^{\text{th}}$ place only, and $0$ everywhere else. We begin by showing that $[\bar\theta, \bar\nu_i]=0$ for each $i=1,\dots, n$. Without loss of generality, take $i=1$, so that $\bar\nu_1=1\times0\times\dots\times 0$. Note that since $0\times0\times\dots\times0\leq \bar\theta\land\bar\nu_1\leq\bar\nu_1=1\times0\times\dots\times0$, Lemma \ref{lemma_product_congruence} guarantees that there  is a $\bar\varphi$ in $\Con(A_1)$ such that $\bar\theta\land\bar\nu_1 = \bar\varphi\times0\times\dots\times 0$. Now since $\bar\theta\restriction_B=0$, we have $(\bar\varphi\times0\times\dots\times0)\restriction_B=0$ as well, and as the embedding $f$ is product--essential, we must have $\bar\varphi =0$. It follows that
$[\bar\theta,\bar\nu_1]\leq \bar\theta\land\bar\nu_1=0$

In the same way, we see that $[\bar\theta, \bar\nu_i]=0$ for all $i=1,\dots, n$. Now
$[\bar\theta,1\times\dots\times 1] = [\bar\theta,\bigvee_{i=1}^n\bar\nu_i] = \bigvee_{i=1}^n[\bar\theta,\bar\nu_i] = 0$.
\eproof

For the next few lemmas we fix the following: We are given a subdirect product--essential embedding $e\colon A\into\prod_{i=1}^n A_i$. For $i=1,\dots, n$, let $\eta_i$ be the kernel of the natural projection $\pi_i\circ e\colon A\onto A_i$, so that $A/\eta_i\cong A_i$. Throughout, we will take $A_i=A/\eta_i$, when convenient. We also assume that the embedding $e\colon A\into\prod_{i=1}^nA_i$ is an inclusion, so that
\[a\in A\qquad\text{is identified with}\qquad (a/\eta_1,\dots, a/\eta_n)\in\prod_{i=1}^n A_i\]Further, define
\[\alpha_i:=\eta_i\lor(\bigwedge_{j\not=i}\eta_j)\in\Con(A)\qquad\qquad\bar\alpha_i:=\alpha_i/\eta_i\in\Con(A_i)\]

\begin{lemma} \label{lemma_alpha_dense} $\bar\alpha_1\times\dots\times\bar\alpha_n$ is dense in $\Con(\prod_{i=1}^nA_i)$.
\end{lemma}

\bproof Let  $\bar\mu_i:=1\times1\dots\times1\times\bar\alpha_i\times1\times\dots\times 1$ in $\Con(\prod_{i=1}^nA_i)$, where $\bar\alpha_i$ occurs in the $i^{\text{th}}$ place, and $1$ everywhere else. Then $\prod_{i=1}^n\bar\alpha_i=\bigwedge_{i=1}^n\bar\mu_i$, so by Lemma \ref{lemma_prod_dense} it suffices to show that $\bar\mu_i$ is dense for each $i=1,\dots,n$.

Without loss of generality, we may assume that $i=1$. Suppose therefore that $\bar\psi\land\bar\mu_1=0$. Then since $0\times1\times\dots\times 1\leq (0\times 1\times\dots\times 1)\lor\bar\psi\leq 1\times1\times\dots\times1$, there is by Lemma \ref{lemma_product_congruence} a $\bar\varphi_1\in\Con(A_1)$ such that $(0\times 1\times\dots\times 1)\lor\bar\psi=\bar\varphi_1\times1\times\dots\times1$. As $A_1=A/\eta_1$, there is $\varphi_1\geq\eta_1$ in $\Con(A)$ such
that $\bar\varphi_1=\varphi_1/\eta_1$. Now 
\[\aligned &\bar\psi\land\bar\mu_1=0\\\Longrightarrow\qquad &(0\times1\times\dots\times 1)\lor(\bar\psi\land\bar\mu_1)=0\times1\times\dots\times1\\\Longrightarrow\qquad &\Big((0\times1\times\dots\times1)\lor\bar\psi\Big)\land\bar\mu_1=0\times1\times\dots\times 1\quad\text{by modularity}\\\Longrightarrow\qquad &(\bar\alpha_1\land\bar\varphi_1)\times1\times\dots\times 1=0\times 1\times\dots\times 1\\	
\Longrightarrow\qquad &\bar\alpha_1\land\bar\varphi_1=0\\\endaligned\]
\[\aligned
\Longrightarrow\qquad &\alpha_1\land\varphi_1\leq\eta_1\\\Longrightarrow\qquad &(\eta_1\lor(\eta_2\land\dots\land\eta_n))\land\varphi_1\leq\eta_1\\
\Longrightarrow\qquad &\varphi_1\land\eta_2\land\dots\land\eta_n\leq\eta_1\quad\qquad\qquad\qquad\qquad\qquad\text{by modularity}\\
\Longrightarrow\qquad &\varphi_1\land\eta_2\land\dots\land\eta_n=0\endaligned\]
Now since $\eta_1\leq\varphi_1$, it follows by Remarks \ref{rem_product_essential}(\ref{rem_eta_max}) that $\varphi_1=\eta_1$. Hence $\bar\varphi_1=0$ in $\Con(A_1)$, i.e. $(0\times1\times\dots\times 1)\lor\bar\psi = 0\times1\times\dots\times 1$. Thus $\bar\psi\leq 0\times1\times\dots\times 1\leq\bar\alpha_1\times 1\times\dots\times 1$, i.e. $\bar\psi\leq\bar\mu_1$, and so $0=\bar\psi\land\bar\mu_1=\bar\psi$.
\eproof

\begin{lemma}\label{lemma_beta_eta} Suppose that for $i=1,\dots, n$ we have $\beta_i\in\Con(A)$ such that 
\begin{enumerate}[(i)]\item $\eta_i\leq \beta_i\leq\alpha_i$, and 
\item $\bar\beta_i:=\beta_i/\eta_i$ is a central congruence in $\Con(A_i)$.
\end{enumerate}
Then for each $k=1,\dots, n-1$ we have \[(\beta_1\land\dots\land\beta_k)\circ\beta_{k+1} = (\eta_1\land\dots\land\eta_k)\circ\eta_{k+1}\]
\end{lemma}

\bproof Since each $\bar\beta_i$ is central in $\Con(A_i)$, the congruence $0\times\dots\times0\times\bar\beta_i\times0\times\dots\times 0$ is central in $\Con(\prod_{i=1}^n A_i)$. As restrictions of central congruences are central, we see that $\xi_i:=(0\times\dots\times0\times\bar\beta_i\times0\times\dots\times 0)\restriction_A$ is central in $\Con(A)$. In particular, each $\xi_i$ is abelian, and hence permutes with every congruence in $\Con(A)$.

Recall that $A_i$ is identified with $A/\eta_i$ and that elements  $a\in A$ are identified with tuples $(a/\eta_1,\dots, a/\eta_n) \in\prod_{i=1}^n A_i$.  Now
\[\aligned &(a,b)\in\xi_i\\\Longleftrightarrow\qquad& (a/\eta_1,\dots,a/\eta_n)\;(0\times\dots\times\bar\beta_i\times\dots\times 0)\;(b/\eta_1,\dots, b/\eta_n)\\\Longleftrightarrow\qquad&(a,b)\in \beta_i\land(\bigwedge_{j\not=i}\eta_j)\endaligned\] and thus
$\xi_i =\beta_i\land(\bigwedge_{j\not=i}\eta_j)$, for all $i=1,\dots, n$.
From this, we obtain the following fact, which we will use several times in the sequel:
\[\xi_i\leq\eta_j\qquad\text{ whenever }\quad i\not=j\]

Next, observe that, since $\eta_i\leq\beta_i\leq\alpha_i$, we have 
\[\beta_i=\beta_i\land \alpha_i=\beta_i\land\Big((\bigwedge_{j\not=i}\eta_j)\lor\eta_i\Big)= (\beta_i\land\bigwedge_{j\not=i}\eta_j)\lor\eta_i=\xi_i\lor\eta_i=\xi_i\circ\eta_i
\]

We will now show by induction that for all $k=1,\dots, n$ we have
\[\beta_1\land\dots\land\beta_k = (\eta_1\land\dots\land\eta_k)\circ\xi_1\circ\dots\circ \xi_k\] We have just shown this is true for $k=1$. Assuming now that $(\beta_1\land\dots\land\beta_{k-1}) = (\eta_1\land\dots\land\eta_{k-1})\circ\xi_1\circ\dots\circ\xi_{k-1}$, we see that \[\aligned \beta_1\land\dots\land\beta_{k-1}\land\beta_k&= [(\eta_1\land\dots\land\eta_{k-1})\lor\xi_1\lor\dots\lor\xi_{k-1}] \land(\eta_k\lor\xi_k)\\&= \Big([(\eta_1\land\dots\land\eta_{k-1})\lor\xi_1\lor\dots\lor\xi_{k-1}]\land\eta_k\Big)\lor\xi_k\\&\text{\qquad\qquad (since $\xi_k\leq\eta_1\land\dots\land\eta_{k-1}$)}\\
&=\Big([\eta_1\land\dots\land\eta_{k-1}\land\eta_k]\lor\xi_1\lor\dots\lor \xi_{k-1}\Big)\lor\xi_k\\&\text{\qquad\qquad (since $\xi_1\lor\dots\lor\xi_{k-1}\leq\eta_{k}$)}\\
&=(\eta_1\land\dots\land\eta_k)\circ\xi_1\circ\dots\circ\xi_k
\endaligned\] using the fact that each $\xi_i$ permutes with every congruence. This completes the induction step.

Now using the permutability of the $\xi_i$ we see that \[\aligned (\beta_1\land\dots\land\beta_k)\circ\beta_{k+1}&= [(\eta_1\land\dots\land\eta_k)\circ\xi_1\circ \dots\circ\xi_k]\circ(\xi_{k+1}\circ\eta_{k+1})\\&=[(\eta_1\land\dots\land\eta_k)\circ\xi_{k+1}]\circ(\xi_1\circ\dots\circ \xi_k\circ\eta_{k+1})\\&=(\eta_1\land\dots\land\eta_k)\circ\eta_{k+1}\endaligned\]
because $\xi_{k+1}\leq\eta_1\land\dots\land\eta_k$ and $\xi_1\lor\dots\lor\xi_k\leq\eta_{k+1}$.
\eproof

\begin{lemma} \label{lemma_[A]theta} Suppose that we have $\beta_i$ as in as in Lemma \ref{lemma_beta_eta}. 
If $\bar\theta\in\Con(\prod_{i=1}^nA_i)$ is such that $\bar\theta\leq\bar\beta_1\times\dots\times\bar\beta_n$, then $[A]\bar\theta=A$.\end{lemma}

\bproof Since every $\bar\theta$--congruence class is contained in a $(\bar\beta_1\times\dots\times\bar\beta_n)$--class, it suffices to show that $[A](\bar\beta_1\times\dots\times\bar\beta_n) = A$. 

Recall once more that in our subdirect product--essential embedding $e\colon A\into\prod_{a=1}^nA_i$,  the algebras $A_i$ are identified with $A/\eta_i$ and  elements  $a\in A$ are identified with tuples $(a/\eta_1,\dots, a/\eta_n) \in\prod_{i=1}^n A_i$. Now suppose  that $(x_1/\eta_1,\dots, x_n/\eta_n)\in[A](\bar\beta_1\times\dots\times\bar\beta_n)$.This means  that there exists an $a\in A$ such that $(x_1/\eta_1,\dots,x_n/\eta_n)\, (\bar\beta_1\times\dots\times\bar\beta_n)\,(a/\eta_1, \dots, a/\eta_n)$, i.e. that $\bigcap_{i=1}^n[x_i]\beta_i\not=\emptyset$. We must show that $(x_1/\eta_1,\dots, x_n/\eta_n)\in A$, i.e. that there exists an element $b\in A$ such that $(x_1/\eta_1,\dots, x_n/\eta_n) = (b/\eta_1,\dots, b/\eta_n)$, i.e. that $\bigcap_{i=1}^n[x_i]\eta_i\not=\emptyset$.

It therefore suffices to prove that if $\bigcap_{i=1}^k[x_i]\beta_i\not=\emptyset$, then $\bigcap_{i=1}^k[x_i]\eta_i\not=\emptyset$, for $k=1,\dots, n$, and we do this by induction on $k$. There is nothing to prove for the case $k=1$. Proceeding with the induction step, suppose that $\bigcap_{i=1}^{k+1}[x_i]\beta_i\not=\emptyset$, and that $a\in 
\bigcap_{i=1}^{k+1}[x_i]\beta_i$. Then also $\bigcap_{i=1}^k[x_i]\beta_{i}\not=\emptyset$, so by induction hypothesis there exists an element $b\in\bigcap_{i=1}^k[x_i]\eta_i$. But then both $a,b\in\bigcap_{i=1}^k[x_i]\beta_i$ so that $b\,(\beta_1\land\dots\land\beta_k)\,a\,\beta_{k+1}\,x_{k+1}$, and hence $b\,(\beta_1\land\dots\land\beta_k)\circ\beta_{k+1}\,x_{k+1}$.
By Lemma \ref{lemma_beta_eta}, we see that there is $c\in A$ such that $b\,(\eta_1\land\dots\land\eta_k)\,c\,\eta_{k+1}\, x_{k+1}$. Then for $i=1,\dots k$ we have $x_i/\eta_i = b/\eta_i = c/\eta_i$. We also have $x_{k+1}/\eta_{k+1} = c/\eta_{k+1}$. Thus $c\in\bigcap_{i=1}^{k+1}[x_i]\eta_i$, so that the latter set is non--empty. This completes the induction step.
\eproof

We can now prove the main result of this section.  The proof adopts a strategy implicit in Davey and Kov\'acs\cite{DK:1986}, to whit that  to prove that an extension $A\into B$ is essential, one may proceed as follows: Given $\alpha\in\Con(B)$ such that $\alpha\restriction_A=0_A$,
\begin{enumerate}[(i)]\item Show that $\alpha$ is central in $\Con(B)$.
\item Then show that there is a congruence $\beta$  which is dense in $\Con(B)$ such that $[A]\beta = A$.
\item Now define $\theta:=\alpha\land\beta$. Then also $\theta\restriction_A=0_A$. Moreover, $\theta\leq\alpha$, so $\theta$ is central. As also $\theta\leq\beta$, we see that $[A]\theta=A$.
\item By Proposition \ref{propn_central_essential}, it follows that $\theta=0_{B}$, i.e. that $\alpha\land\beta =0_B$. As $\beta$ is dense, it follows that $\alpha=0_B$, and hence  $A\into B$ is an essential extension.
\end{enumerate}

\noindent{\em Proof of Theorem \ref{thm_sd_pr_ess}:} Suppose that $\bar\theta\in\Con(\prod_{i=1}^n A_i)$ is such that $\bar\theta\restriction_A=0$.  Define $\bar\psi:=\bar\theta\land(\bar\alpha_1\times\dots\times\bar\alpha_n)$, so that also $\bar\psi\restriction_A=0$. By Proposition \ref{propn_theta_central} it follows that $\bar\psi$ is a central congruence in $\Con(\prod_{i=1}^nA_i)$, and thus that $\bar\psi\leq\bar\beta_1\times\dots\times\bar\beta_n$, where $\bar\beta_i:=\bar\alpha_i\land\zeta_{A_i}$ and  $\zeta_{A_i}$ denotes the center of $A_i$. Hence by Lemma \ref{lemma_[A]theta} we may conclude that $[A]\bar\psi = A$. It now follows from Proposition \ref{propn_central_essential} that $\bar\psi = 0$, i.e. that $\bar\theta\land(\bar\alpha_1\times\dots\times\bar\alpha_n) = 0$. Finally, Lemma \ref{lemma_alpha_dense} shows that $\bar\theta=0$, as required. \hfill $\Box$

\section{Absolute Retracts}
With Theorem \ref{thm_sd_pr_ess} in hand, we can investigate absolute retracts in congruence modular varieties. To start, we obtain the following result:
\begin{theorem}\label{thm_DK} {\rm (Davey and Kov\'acs\cite{DK:1986})} If $A$ is a directly indecomposable absolute retract in a congruence modular variety, then it is finitely subdirectly irreducible. 
\end{theorem}
\bproof Suppose that $A$ is an absolute retract. If $A$ is not finitely subdirectly irreducible, then there exist $\eta_1,\eta_2>0$ in $\Con(A)$ such that $\eta_1\land\eta_2=0$. As in Remarks \ref{rem_product_essential}(\ref{rem_sd_pr_ess}), we may assume that $\eta_1,\eta_2$ are maximal with respect to the property of having zero meet, so that the canonical embedding $A\into A/\eta_1\times A/\eta_2$ is product essential, and hence essential. But an absolute retract has no proper essential extensions, so $A\cong A/\eta_1\times A/\eta_2$, so that $A$ is directly decomposable --- contradiction. Hence $A$ is finitely subdirectly irreducible.
 \eproof

\begin{theorem}\label{thm_dec_AR} If $A$ is an absolute retract in a congruence modular variety, and if $\Con(A)$ is finite dimensional, then $A$ is a finite product of subdirectly irreducible algebras.
\end{theorem}

\bproof  Suppose that $0_A=\bigwedge_{i=1}^n\eta_i$ is a representation of $0_A$ as an irredundant meet (where the $\eta_i$ are not assumed to be meet irreducible). Then $\eta_1>\eta_1\land\eta_2>\dots>\bigwedge_{i=1}^n\eta_i$, and hence $n\leq\text{height(Con($A$))}$. Now let $m$ be the maximum integer for which there exists an irredundant meet representation $0_A=\bigwedge_{i=1}^m\eta_i$ of length $m$. By the procedure outlined in Remarks \ref{rem_product_essential}(\ref{rem_sd_pr_ess}), we may choose $\varphi_i\geq \eta_i$ (for $1\leq i\leq m$) so that \begin{enumerate}[(i)]\item $\bigwedge_{i=1}^m\varphi_i=0_A$; this meet is then clearly irredundant.
\item If $\theta_i\geq\varphi_i$ for $1\leq i\leq m$ are such that $\bigwedge_{i=1}^m\theta_i=0_A$ then each $
\theta_i=\varphi_i$.
\end{enumerate} It then follows that the natural subdirect embedding $A\into\prod_{i=1}^nA/\varphi_i$ is product--essential, and hence essential. But then $A\cong \prod_{i=1}^n A/\varphi_i$, as $A$ does not have any proper essential extensions. 

Moreover, each $\varphi_i$ is meet irreducible. To see this, suppose for example that $\varphi_1=\chi\land\psi$, where $\chi,\psi>\varphi_1$.  Then $0_A=\chi\land\psi\land\bigwedge_{i=2}^m\varphi_i$. Since $\chi\land\bigwedge_{i=2}^m\varphi_i>0_A$ and $\psi\land\bigwedge_{i=2}^m\varphi_i>0_A$, we see that the representation $0_A=\chi\land\psi\land\bigwedge_{i=2}^m\varphi_i$ is an irredundant meet of length $m+1$, contradicting the maximality of $m$. Hence $\varphi_1$ is meet irreducible, and as $\Con(A)$ has finite height, it is completely meet irreducible. By the same argument, it follows that each $\varphi_i$ is completely meet irreducible, so that $A\cong\prod_{i=1}^m A/\varphi_i$ decomposes $A$ into a product of subdirectly irreducible algebras. 
\eproof

Davey and Kov\'acs\cite{DK:1986} actually prove a stronger result than Theorem \ref{thm_DK}, namely that a directly indecomposable absolute retract in a congruence modular variety is a finitely subdirectly irreducible which is either centerless or abelian.  We will now partly extend this result to include products, but with the added assumption that the base variety is residually small. Now any residually small congruence modular variety satisfies the commutator identity (C1): $\alpha\land[\beta,\beta]=[\alpha\land\beta, \beta]$. Observe that, since $\zeta_{A\times B} = \zeta_A\times\zeta_B$ and $[1_A\times 1_B, 1_A\times 1_B] = [1_A,1_A]\times [1_B,1_B]$, the properties of being {\em centerless} and of being {\em abelian} are both preserved under finite products.

\begin{theorem}\label{thm_AR_prod_ab_cless} Let $A$ be an absolute retract in a congruence modular variety which satisfies (C1). Then $A$ is a product of a centerless algebra and an abelian algebra. In particular, if $A$ is directly indecomposable, then $A$ is either centerless or abelian.
\end{theorem}
\bproof
Note that $0=[\zeta,1]=\zeta\land[1,1]$. Now choose $\theta\geq\zeta$ and $\psi\geq [1,1]$ in $\Con(A)$ maximal so that $\theta\land\psi =0$, as discussed in Remarks \ref{rem_product_essential}(\ref{rem_sd_pr_ess}). Then  the subdirect embedding $A\into A/\theta\times A/\psi$ is product--essential, hence essential. Since $A$ is an absolute retract, we have $A\cong A/\theta\times A/\psi$. But $[\theta,1]= \theta\land[1,1]=0$, so $\theta$ is central, i.e. $\theta =\zeta$. 

Furthermore, $A/\zeta$ is centerless: For if $\xi\geq\zeta$ in $\Con(A)$ is such that $\xi/\zeta$ is central in $\Con(A/\zeta)$, then $[\xi/\zeta, 1/\zeta]= 0$, so $[\xi,1]\leq\zeta$. But then $0=[[\xi,1],1] =[\xi,1]\land [1,1]=[\xi,1]$, and hence $\xi$ is central, so $\xi=\zeta$. Thus  $\xi/\zeta=0$ in $\Con(A/\zeta)$.

Finally, since $\psi\geq[1,1]$, $A/\psi$ is abelian. Hence $A\cong A/\zeta\times A/\psi$ is the product of a centerless and an abelian algebra.
\eproof

\begin{theorem} Suppose that $A$ is an absolute retract in a congruence modular variety satisfying (C1), such that $\Con(A)$ is finite dimensional. If $A$ satisfies the unique factorization property, then $A$ is a finite product of subdirectly irreducible algebras, each of which is either centerless or abelian. \newline In particular, this conclusion is valid when either (i) $A$ is congruence --permutable with a one--element subalgebra, or  when (ii) $\Con(A)$ is finite.\end{theorem}

\bproof By Theorem \ref{thm_dec_AR}, $A$ is a product of subdirectly irreducible algebras $A=\prod_{i=1}^kA_i$. Furthermore, $A=A/\zeta\times A/\psi$ is also a product of a centerless and an abelian algebra, by Theorem \ref{thm_AR_prod_ab_cless}. Now subdirectly irreducible algebras are directly indecomposable, and hence, by unique factorization,  $A/\zeta$ must be the product of some of the $A_i$, and $A/\psi$ the product of the remaining $A_i$. By reindexing, we may assume that $A/\zeta=\prod_{i=1}^n A_i$, and that $A/\psi =\prod_{i=n+1}^k A_i$.  Now since $A/\zeta$ is centerless and $A/\psi$ is abelian, we see that $A_i$ is centerless when $i\leq n$, and abelian when $n+1\leq i\leq k$.

The Birkhoff--Ore Theorem states that if $A$ has a one--element subalgebra, and a finite--dimensional congruence lattice with permuting congruences, then $A$ has the unique factorization property. A theorem of J\'onsson reaches the same conclusion in the case that $\Con(A)$ is modular and finite. See Chapter 5 of \cite{MMT:1987} for a proof of both assertions.\eproof

\begin{remarks}\rm
\begin{enumerate}[(a)]\item
The above--mentioned results of Birkhoff--Ore and J\'onsson admit a common generalization: {\em Is every algebra with finite dimensional modular congruence lattice is uniquely factorable?} This problem is still open, cf. \cite{MMT:1987}.
\item It is easy to show that if $A=\prod_{i\in I}A_i$ is an absolute retract with a one--element subalgebra, then each $A_i$ is an absolute retract also. In that case we can deduce from Theorems \ref{thm_AR_prod_ab_cless}  and \ref{thm_dec_AR} that in a residually small congruence modular variety $\mathcal V$:
\begin{enumerate}[(i)]\item If $A$ is an absolute retract, then $A/\zeta$ is an absolute retract. Moreover, if $\zeta\not=0$, then some non--trivial abelian image of $A$ is an absolute retract.
\item If $A$ is an absolute retract with finite--dimensional congruence lattice in a congruence modular variety, then $A$ is a finite product of  {\em maximal} subdirectly irreducibles, each of which is either centerless or abelian.
\end{enumerate}
\item Note that every (weakly) injective algebra is an absolute retract. Furthermore, a finite algebra is algebraically closed in a variety if and only if it is an absolute retract (cf. \cite{Ouwehand:2009}). Thus the above results also have implications for (weakly) injective and algebraically closed algebras. \end{enumerate}
\end{remarks}

\appendix\section{Auxiliary Results, Not For Publication}
The following result is from Davey and Kov\'acs\cite{DK:1986}, who base their argument on a quite complicated isomorphism (obtained by Gumm\cite{Gumm:1983}) between the lattice of central congruences of an algebra $A$ and the lattice of subalgebras of a certain algebra defined on a congruence class of the center of $A$. The direct proof given here incorporates arguments of from Gumm\cite{Gumm:1983}.
\begin{proposition}\label{propn_central_extend} Suppose that $\bar A$ is an algebra in a congruence modular variety, and that $A$ is a subalgebra of $\bar A$. 
Let $a\in A$, and let 
$\bar\alpha\leq\bar \zeta$ in $\Con(\bar A)$, where $\bar \zeta$ is the center of $\bar A$. Suppose further that $\beta\leq\bar\alpha\restriction_A$ in $\Con(A)$. Now define
\[\bar\beta:=\{(x,y)\in\bar\zeta: d(a,x,y)\in [a]\beta\}\]
Then $\bar\beta\in\Con(\bar A)$, $\bar\beta\leq\bar\alpha$ and $\bar\beta\restriction_A = \beta$.
\end{proposition}

\bproof
Suppose that $(x,y)\in\bar\beta$. Then $d(a,x,y)\in[a]\beta\subseteq[a]\bar\alpha$ and hence \\$x = d(x,a,a) \,\bar\alpha\, d(x,a,d(a,x,y)) = y$. But $d(x,a,d(a,x,y)) = y$ by Corollary \ref{corr_p_id}, as $\bar\beta\subseteq\bar\zeta$. Thus also $(x,y)\in\bar\alpha$, which shows that $\bar\beta\subseteq\bar\alpha$.

Let $(x,y)\in\beta$. Then $d(a,x,y)\,\beta, d(a,x,x) = a$, and hence $(x,y)\in\bar\beta$, i.e. $\beta\subseteq\bar\beta\restriction_A$. Conversely, if $(x,y)\in\bar\beta\restriction_A$, then $d(a,x,y)\,\beta, a$ and hence $x=d(x,a,a)\,\beta\,d(x,a,d(a,x,y)) =y$ by Corollary \ref{corr_p_id}. Thus also $\bar\beta\restriction_A\subseteq \beta$.

It remains to show that $\bar\beta$ is a congruence on $\bar A$. As $d(a,x,x) = a$, it is clear that $\bar\beta$ is reflexive. 

Now suppose that $(x,y)\in\bar\beta$, i.e. that $d(a,x,y)\in[a]\beta\subseteq A$. Since $\bar\beta\leq\bar\zeta$, we see that $d$ commutes with itself on $\left(\begin{matrix}a&a&a\\a&x&y\\a&x&x\end{matrix}\right)$, and conclude that $d(a,y,x)=d(a,d(a,x,y),a)$.
Now as $a,d(a,x,y)\in A$, it follows that $d(a,y,x)\in A$ also. Clearly $d(a,d(a,x,y),a)\,\beta\, d(a,a,a)$ yields $d(a,y,x)\in[a]\beta$. Hence $\bar\beta$ is symmetric.

Next, suppose that $(x,y), (y,z)\in\bar\beta$. Since $\bar\beta\leq\bar\zeta$, $d$ commutes with itself on $\left(\begin{matrix}a&x&y\\a&y&y\\a&y&z\end{matrix}\right)$ and hence $d(a,x,z)= 
d(d(a,x,y),a,d(a,y,z))$. Since $d(a,x,y),d(a,y,z)\in [a]\beta\subseteq A$, see that $d(a,x,z)\in A$ also. Then $d(d(a,x,y),a,d(a,y,z))\,\beta\, d(a,a,a)$ yields that $d(a,x,z)\in[a]\beta$. Hence $\bar\beta$ is transitive.

That $\bar\beta$ is compatible with all the operations follows once again from Theorem \ref{thm_Gumm_9.1}: Suppose that $f$ is an $n$--ary operation, and that $(x_i,y_i)\in\bar\beta$ for $i=1,\dots, n$. Let $\mathbf x:=(x_1,\dots, x_n)$, with $\mathbf y$ defined similarly, and let $\mathbf a = (a,a,\dots, a)$ denote an $n$--tuple consisting entirely of $a$'s. Since $\bar\beta\leq\bar\zeta$, we have $f(\mathbf x)\,\bar\zeta\,f(\mathbf y)$, and hence $d$ commutes with itself on $\left(\begin{matrix} a&f(\mathbf x)&f(\mathbf x)\\
f(\mathbf a)&f(\mathbf x)&f(\mathbf x)\\
f(\mathbf a)&f(\mathbf x)&f(\mathbf y)\end{matrix}\right)$, so that $d(a,f(\mathbf x), f(\mathbf y)) =d(a, f(\mathbf a), d(a,f(\mathbf x),\mathbf f(y))$. But as $f$ commutes with $d$ on $\left(\begin{matrix}a&x_1&y_1\\\vdots&\vdots&\vdots\\a&x_n&y_n\end{matrix}\right)$, we have $d(a,f(\mathbf x),f(\mathbf y) )= f(d(a,x_1,y_1),\dots d(a,x_n,y_n))$, so that 
$d(a,f(\mathbf x),f(\mathbf y) = d\Big(a, f(\mathbf a), f(d(a,x_1,y_1),\dots, d(a,x_n,y_n))\Big)$.
Now as $d(a,x_i,y_i)\in[a]\beta$ for $i=1,\dots, n$, it follows that $d(a,f(\mathbf x), f(\mathbf y))\,\beta \,d(a, f(\mathbf a), f(\mathbf a)) =a$, and thus that $f(\mathbf x)\,\bar\beta\,f(\mathbf y)$. Hence $\bar\beta\in\Con(\bar A)$.
\eproof

%XXXXXXXXXXXXXXXXXXXXXXXXXXXXXXXXXXXXXXXXXXXXXXXXXXXXXXXXXXXXXXXXXXXXXXX

Our next proposition is again from Davey and Kov\'acs\cite{DK:1986}, with only slight modifications. We include the proof to keep this exposition self--contained.

\begin{proposition}\label{propn_central_cube} Suppose that $A$ is non--abelian, and that its center $\zeta$ is dense in $\Con(A)$. Then $A$ has a proper essential extension. 
\end{proposition}
\bproof
Clearly $B:=\{(a,b,c)\in A^3:a\,\zeta\,b\,\zeta\,c\}$ is a subalgebra of $A$. By Theorem \ref{thm_Gumm_9.1}, $d:B\onto A$ is a surjective homomorphism. Let $\theta:=\ker d$.  By Proposition \ref{propn_central_extend}, there is a central congruence $\Theta\leq\zeta^3$ such that $\Theta\restriction_B=\theta$. We use the Davey and Kov\'acs \cite{DK:1986} strategy decribed just prior to the proof of Theorem \ref{thm_sd_pr_ess} to prove that the induced embedding $B/\theta\into A^3/\Theta$ is an essential extension.

(i) Suppose that $\Psi\geq\Theta$ in $\Con(A^3)$ satisfies $\Psi/\Theta\land\zeta^3/\Theta=0$, i.e. that $\Psi\land\zeta^3=\Theta$.  We begin by showing that $\Psi$ is a central congruence. By Lemma \ref{lemma_product_congruence} there is a $\beta\in\Con(A)$ such that $\Psi\land (0\times 1\times 0) = 0\times\beta\times0$, so that $0\times(\beta\land\zeta)\times 0 \leq \Psi\land\zeta^3=\Theta$. But if $(a,b)\in \beta\land\zeta$, then $(a,a,a,)\,0\times(\beta\land\zeta)\times 0\,(a,b,a)$, so that $(a,a,a)\,\Theta\,(a,b,a)$. But as $(a,a,a), (a,b,a)\in B$, it follows that $(a,a,a)\,\theta\,(a,b,a)$, and hence $b=d(b,a,a)=d(b,a,d(a,b,a))=a$, by Corollary \ref{corr_p_id}. It follows that $\beta\land\zeta=0$, and thus that $\beta=0$, as $\zeta$ is dense in $\Con(A)$. Thus $\Psi\land(0\times 1\times 0) = 0$, and hence $[\Psi, 0\times 1\times 0] = 0$. In an analogous --- but slightly simpler --- fashion does it follow that $[\Psi, 1\times 0\times 0]=0=[\Psi, 0\times0\times 1]$. Hence $[\Psi, 1\times 1\times 1]=0$, so that $\Psi$ is a central congruence on $A^3$.

(ii) It follows that $\Psi\leq \zeta^3$, and hence that $\Theta=\Psi\land\zeta^3=\Psi$. We thus see that if $\Psi/\Theta\land \zeta^3/\Theta=0$ in $\Con(A^3/\Theta)$, then $\Psi/\Theta=0$. Hence $\zeta^3/\Theta$ is dense in $\Con(A^3/\Theta)$. Moreover, we clearly have $[B]\zeta^3 =B$, and thus also $[B/\theta]\zeta^3/\Theta= B/\theta$. 

With (i) and (ii) of the Davey and Kov\'acs strategy satisfied, (iii) and (iv) follow. Hence $B/\theta\into A^3/\Theta$ is an essential extension. As $\Theta\leq\zeta^3$, this embedding is surjective  exactly when $B=A^3$, i.e. when $\zeta=1$. But then $A$ is abelian, contradiction.\eproof

It now follows easily that a directly indecomposable absolute retract $A$ is either centerless or abelian. For if $A$ is not centerless, then as $A$ is finitely subdirectly irreducible, we see that $\zeta_A$ is dense in $\Con(A)$. But then as $A$ has no proper essential extension, it must be abelian.
\bibliographystyle{plain}
\bibliography{bib_Modular_AR}

\end{document}